\documentclass{amsart}
\begin{document}
\title{Warped products and Reissner-Nordstrom metric}
\author{Soon-Tae Hong$^{1}$, Jaedong Choi$^{2}$ and Young-Jai Park$^{1}$}
\address{$^1$Department of Physics and Basic Science Research
  institute, Sogang University\\
 C.P.O. Box 1142, Seoul 100-611, Korea }
\email{sthong@sogang.ac.kr\\ yjpark@sogang.ac.kr}
\address{$^2$Department of Mathematics, P.O. Box 335-2, Air Force Academy, Ssangsu, Namil\\
Cheongwon, Chungbuk 363-849, Korea}\email{\noindent: jdong@afa.ac.kr}

\noindent
\begin{abstract} 
  
We study a multiply warped products manifold associated with the Reissner-Nordstrom 
metric to investigate the physical properties inside the black hole 
event horizons.  It is shown that, different from the uncharged 
Schwarzschild metric, the Ricci curvature components inside the 
Reissner-Nordstrom black hole horizons are nonvanishing, while the 
Einstein scalar curvature vanishes even in the interior of the charged metric.  
Introducing a perfect fluid inside the Reissor-Nordstrom black hole, it is also shown that 
the charge plays effective roles of decreasing the mass-energy density and the pressure of the 
fluid inside the black hole. 

\end{abstract} \vskip10pt  

\maketitle
\noindent
{\bf I. Introduction}
\vskip10pt

Since the pioneering work in 1976${}^1$, thermal
Hawking effects on a curved manifold${}^2$ have been studied 
as an Unruh effect in a higher flat dimensional spacetime.
Following the global embedding Minkowski space approach${}^{3-6}$, several 
authors recently have shown that this approach could yield a unified
derivation of temperature for various curved manifolds in (2+1) dimensions${}^{7-13}$ 
and in (3+1) dimensions${}^{7,14-16}$.  However all these higher 
dimensional embedding solutions have been constructed outside the event 
horizons of the metrics.  
\vskip8pt

On the other hand, the concept of a warped products manifold was 
introduced by Bishop and O'Neill${}^{17}$, where it served to provide a class of complete 
Riemannian manifolds with everywhere negative curvature${}^{18}$. 
The connection with general relativity was first made by Beem, Ehrich, and Powell, who pointed 
out that several of the well-known exact solutions to Einstein field equations are 
pseudo-Riemannian warped products$^{19}$. Beem and Ehrich further explored the extent to which certain 
causal and completeness properties${}^{20}$ of a spacetime maybe determined by the presence of 
a warped products structure. After first developing the general theory of warped products to spaces, 
O'Neill then applied the theory to discuss, in turn, the special cases of Robertson-Walker and 
Schwarzshild spacetime. The role of warped products in the study of exact solutions to 
Einstein's equations is now firmly established, and it appears that these structures are generating 
interest in other areas of geometry.  Very recently, Choi has represented the interior Schwarzschild 
spacetime as a multiply warped products spacetime with warping functions$^{21}$ to yield the 
{\it Ricci curvature} in terms of $f_1$, $f_2$ for the multiply warped products of the form 
$M=R^1 \times_{f_1}R^1\times_{f_2} S^2$. 
\vskip8pt

In this paper we will further analyze the multiply warped products manifold associated 
with the Schwarzschild metric by including the charge degrees of freedom so 
that we can investigate the physical properties inside the event 
horizons of the Reissner-Nordstrom black hole.  
\vskip8pt

In section II, we will briefly recapitulate the 
multiply warped products manifold and the 
corresponding Ricci curvature components. In section III, we will 
apply the multiply warped products manifold scheme to the Reissner-Nordstrom metric 
to explicitly obtain the Ricci and Einstein curvatures inside the event 
horizons of the metric and to study the corresponding charge effects.  Introducing a 
perfect fluid inside the Reissor-Nordstrom black hole, we will discuss the charge 
effects 
on the mass-energy density and the pressure of the fluid inside the black hole.

\vskip20pt

\noindent
{\bf II.  Multiply warped products and Ricci curvature}\vskip10pt

In this section, we briefly review a multiply warped products manifold to investigate Ricci 
curvature inside the black hole horizons.
\vskip8pt

A Lorentzian manifold $(M,g)$ is a connected smooth manifold of $d$-dimension $(d\geq 2)$ with a countable basis 
together with a Lorentzian metric \ $g$ \ of signature $(-, +, +,...,+)$. Let $(F_i, g_i)$ be Riemannian 
manifolds, and let $(B, g_B)$ be either a spacetime, or let $B$ be $R^1$ with $g_B=-dt^2$. Let $f_i>0$, 
$i=1,...,n$ be smooth functions on $B$. A {\it multiply warped products spacetime} with base $(B, g_B)$, 
fibers $(F_i,g_i)$ $i=1,...,n$ and warping functions $f_i>0$ is the product manifold 
$(B\times F_1\times...\times F_n, g)$ endowed with the Lorentzian  metric:
$$g=\pi_B^{\ast}g_B+\sum_{i=1}^n(f_i\circ\pi_B)^2\pi_i^\ast g_i\equiv-dt^2+\sum_{i=1}^nf_i^2g_i
\eqno(2.1)
$$
where $\pi_B$, $\pi_i$\ $i=1,...,n$ are the natural projections of $B\times F_1\times...\times F_n$ onto 
$B$ and $F_1$,...,$F_n$, respectively.\vskip8pt 

Thus, warped product spaces are extended to richer class of spaces involving multiply products. Multiply 
warped products spaces were studied by Flores and S\'{a}nchez$^{22}$.   The conditions of spacelike boundaries 
in the multiply warped products spacetimes were studied by Harris${}^{23}$.  The Kasner metric was 
studied as a cosmological model by Sch\"{u}cking and Heckmann${}^{24}$.  Choi has investigated the curvature 
of a multiply warped product with $C^0$-warping functions$^{21}$.
\vskip8pt

From a physical point of view, these spacetimes are interesting, first, because they include classical 
examples of spacetimes: when $n=1$ they are {\it generalized Robertson-Walker} spacetimes, standard 
models of cosmology; when $n=2$ the intermediate zone of Reissner-Norsdstr\"{o}m spacetime and 
interior of Schwarzschild spacetime appear as particular cases.  \vskip20pt

\noindent
{\bf III.  Reissner-Nordstrom black hole as a multiply warped product manifold}
\vskip10pt

In this section, to investigate a multiply warped product manifold for the Reissner-Nordstrom interior solution, 
we start with the four-metric inside the horizon   
$$
ds^{2}=N^{2}dt^{2}-N^{-2}dr^{2}+r^{2}d\Omega^{2}
\eqno(3.1)
$$
where the lapse function for the interior solution is given by
$$N^{2}=-1+\frac{2m}{r}-\frac{Q^{2}}{r^{2}},\eqno(3.2)$$
with a mass $m$ and a charge $Q$, and $d\Omega^{2}=d\theta^{2}+\sin^{2}\theta d\phi^{2}$.  Note that, 
for the nonextremal case, there exist two event horizons $r_{\pm}(Q)$ satisfying the equations 
$0=-1+2m/r_{\pm}-Q^{2}/r_{\pm}^{2}$ such that
$$r_{\pm}=m\pm(m^{2}-Q^{2})^{1/2}.\eqno(3.3)$$
Furthermore the lapse function can be rewritten in terms of these outer and inner horizons 
as follows
$$N^{2}=\frac{(r_{+}-r)(r-r_{-})}{r^{2}}\eqno(3.4)$$
which, for the interior solution, is well defined in the region $r_{-}<r<r_{+}$.  
\vskip20pt 

{\bf Proposition 3.1} \ Let $M$ be a manifold with the Reissor-Nordstrom metric solution 
$ds^{2}=N^{2}dt^{2}-N^{-2}dr^{2}+r^{2}d\Omega^{2}$ (where $d\Omega^2=d\theta^2+{\sin^2{\theta}}
\thinspace d\phi^2)$, then we have $M$ as a multiply warped manifolds  with warping functions $f_1,\ f_2$ 
$$ds^{2}=-d\mu^{2}+f_{1}^{2}(\mu)d\nu^{2}+f_{2}^{2}(\mu)d\Omega^{2}
\eqno(3.5)$$
where 

$$\mu=2m\cos^{-1}\left(\frac{r_{+}-r}{r_{+}-r_{-}}\right)-(r_{+}-r)^{1/2}(r-r_{-})^{1/2}=F(r),
\eqno(3.6)$$
$$f_{1}(\mu)=\ \left(-1+\frac{2m}{F^{-1}(\mu)}-\frac{Q^{2}}{F^{-2}(\mu)}\right)^{1/2},\eqno(3.7)$$
$$f_{2}(\mu)=\ F^{-1}(\mu).\eqno(3.8)
$$
\vskip10pt

{\it Proof:}  Define a new coordinate $\mu$ by
$$d\mu^{2}=N^{-2}dr^{2},\eqno(3.9)$$
which can be integrated to yield
$$\mu=\int_{r_{-}}^{r}\frac{dx~x}{(r_{+}-x)^{1/2}(x-r_{-})^{1/2}}+\mu(r_{-}).\eqno(3.10)$$
We establish $\mu(r_{-})=0$ to yield the analytic solution (3.6).  Note that in (3.6) 
$dr/d\mu >0$ implies $F^{-1}(\mu)$ is well-defined function.  Exploiting the above new coordinate 
(3.6) and redefinition $\nu=t$, we rewrite the metric (3.1) as a multiply warped product $M$ as 
in the form of (2.1) to yield (3.5) with $f_{1}$ and $f_{2}$ in (3.7) and (3.8).
\qed
\vskip20pt

Note that in the vanishing $Q$ limit the above solution is reduced to that of the uncharged 
Schwarzschild case$^{21}$.  Moreover, we have the following boundary conditions for $F(r)$

$${\rm lim}_{r\rightarrow r_{+}}F(r)=m\pi,\hskip0.7cm 
{\rm lim}_{r\rightarrow r_{-}}F(r)=0.
\eqno(3.11)$$

After some algebra, we obtain the following nonvanishing Ricci curvature components
\begin{equation}
\begin{split}\hskip2cm
&R_{\mu\mu}=-\frac{f_{1}^{''}}{f_{1}}-\frac{2f_{2}^{''}}{f_{2}},
\nonumber\\
&R_{\nu\nu}=\frac{2f_{1}f_{1}^{'}f_{2}^{'}}{f_{2}}+f_{1}f_{1}^{''},
\nonumber\\
&R_{\theta\theta}=\frac{f_{1}^{'}f_{2}f_{2}^{'}}{f_{1}}+f_{2}f_{2}^{''}
+f_{2}^{' 2}+1,
\nonumber\\
&R_{\phi\phi}=\left(\frac{f_{1}^{'}f_{2}f_{2}^{'}}{f_{1}}+f_{2}f_{2}^{''}
+f_{2}^{'2}+1 \right)\sin^{2}\theta,\hskip3cm(3.12)\end{split}\nonumber\end{equation}
which hold also in the Schwarzschild metric case$^{21}$.  Note that, as shown in (3.6)-(3.8), 
the argument $\mu$ of $f_{1}$ and $f_{2}$ in (3.12) and $f_{1}$ itself are however 
different from those of the Schwarzschild black hole, since $\mu$ is described in terms of 
$r_{\pm}$ possessing the charge and $f_{1}$ has the additional charge term.   
\vskip20pt

{\bf Proposition 3.2} \ Let $M$ be a multiply warped product manifold 
$M=R^1 \times_{f_1}R^1\times_{f_2} S^2$ for the Reissor-Nordstrom metric solution 
$ds^{2}=-d\mu^{2}+f_{1}^{2}(\mu)d\nu^{2}+f_{2}^{2}(\mu)d\Omega^{2}$ (where 
$d\Omega^2=d\theta^2+{\sin^2{\theta}}
\thinspace d\phi^2)$ with warping functions $f_1,\ f_2$, then we have Ricci curvature components 
as follows,
\begin{equation}
\begin{split}
\hskip1cm&(1)\ R_{\mu\mu}=\frac{Q^{2}}{f_{2}^{4}},\hskip1cm
(2)\ R_{\nu\nu}= -\frac{Q^{2}f_{1}^{2}}{f_{2}^{4}},\\
&(3)\ R_{\theta\theta}
= \frac{Q^{2}}{f_{2}^{2}},\hskip1cm  
(4)\ R_{\phi\phi}=\ \frac{Q^{2}}{f_{2}^{2}}\sin^{2}\theta.
\end{split}\hskip4.0cm (3.13)\nonumber
\end{equation}\vskip10pt

{\it Proof:} Using the explicit expressions for $f_{1}$ and $f_{2}$ in (3.7) and (3.8), we can obtain 
identities for $f_{1}$, $f_{1}^{'}$ and $f_{1}^{''}$ in terms of $f_{1}$, $f_{2}$ 
and their derivatives
\begin{equation}
\begin{split}
\hskip2cm&f_{1}=\ f_{2}^{'},\nonumber\\
&f_{1}^{'}=\ -\frac{m}{f_{2}^{2}}+\frac{Q^{2}}{f_{2}^{3}},\nonumber\\
&f_{1}^{''}=\ -\frac{2f_{1}f_{1}^{'}}{f_{2}}-\frac{Q^{2}f_{1}}{f_{2}^{4}}.
\hskip6cm(3.14)\end{split}\nonumber\end{equation}
Substituting (3.14) into (3.12), we evaluate $R_{\mu\mu}$ 
as follows,
$$
R_{\mu\mu}= 2\big(\frac{f_{1}^{'}}{f_{2}}+\frac{Q^{2}}{2f_{2}^{4}}\big)-\frac{2f_{2}^{''}}{f_{2}}
= \frac{Q^{2}}{f_{2}^{4}}.
$$
Similarly, we establish the other components $R_{\nu\nu}$, $R_{\theta\theta}$ 
and $R_{\phi\phi}$. 
\qed
\vskip20pt

Here one notes that, differently from the Schwarzschild case where we have 
the flat Ricci curvature components as shown with $Q=0$, we have the nonvanishing Ricci tensor 
components for $r_{-}<r<r_{+}$.  Moreover, it is amusing to see that, exploiting the 
metric (3.5), the Einstein scalar curvature is given as follows 
$$R=0,\eqno(3.15)$$
even in the interior of the charged Reissner-Nordstrom black hole horizons.  However, due to the charge 
of the black hole, the vacuum solution in the uncharged Schwarzschild case does not hold any more as 
shown in (3.13).

\vskip8pt 
 
Next, in order to investigate the charge effects inside the Reissor-Nordstrom black hole, 
we assume a perfect fluid which is a continuous distribution of matter with 
stress-energy tensor $T_{ab}$ of the form,
$$
T_{ab}=\rho u_{a}u_{b}+P(g_{ab}+u_{a}u_{b}),
\eqno(3.16)$$
where $u^{a}$ is a unit timelike 4-velocity vector field of the fluid, 
$\rho$ is the mass-energy density and $P$ is the pressure of the fluid 
as measured in its rest frame.
\vskip20pt
  
{\bf Proposition 3.3} \ In the interior ($r_{-}<r<r_{+}$) of charged Reissor-Nordstrom 
metric solution, we obtain the following equations.
\begin{equation}
\begin{split}
\hskip1cm&(1)\ \ R_{\mu\mu}-8\pi T_{\mu\mu}=\frac{Q^{2}}{f_{2}^{4}}-8\pi Pf_{1}^2=0,\ \ \ \ \
\\
&(2)\ \ R_{\nu\nu}-8\pi T_{\nu\nu}=-\frac{Q^{2}f_{1}^{2}}{f_{2}^{4}}+8\pi\rho=0,\\
&(3)\ \ R_{\theta\theta}-8\pi T_{\theta\theta}=\frac{Q^{2}}{f_{2}^{2}}-8\pi Pf_{2}^2=0,\\ \ \ \ \
&(4)\ \ R_{\phi\phi}-8\pi T_{\phi\phi}=\frac{Q^{2}}{f_{2}^{2}}\sin^{2}\theta 
-8\pi Pf_{2}^2 \sin^{2}\theta=0.
\end{split}
\hskip2.9cm 
(3.17)
\nonumber\end{equation}
\vskip10pt

{\it Proof:} Substituting into the Einstein equations of motion
$$R_{ab}-\frac{1}{2}Rg_{ab}=8\pi T_{ab},
\eqno(3.18)$$
the Ricci components given by (3.13), the Einstein scalar curvature (3.15) and 
the stress-energy tensor (3.16) given in terms of the warped products, we obtain the above 
results.\qed
\vskip20pt

Here one notes that in the case of the vanishing $Q$ limit which is the Schwarzschild case, one 
has contributions from the mass-energy 
density and the pressure of the fluid in the interior of the black hole.  On the other hand, 
in the charged  Reissor-Nordstrom case, the existence of nonvanishing charge $Q$ plays effective 
roles of 
decreasing the mass-energy density and the pressure of the fluid inside the black hole, 
since the terms associated with the charge are positive definite in the Einstein equations 
of motion.  From a physical viewpoint, the above result is quite consistent with 
the phenomenology that the charge generates a repulsive force inside the charged black hole 
to dilute the density in the perfect fluid itself.      

\vskip20pt

\noindent
{\bf IV.  Conclusions}\vskip10pt
\vskip5pt

We have studied  a multiply warped product manifold associated with the Reissor-Nordstrom metric 
to evaluate the Ricci curvature components 
inside the charged black hole horizons.  Differently from the 
uncharged Schwarzschild metric where 
both the Ricci and Einstein curvatures vanish inside the horizon, the Ricci 
curvature components inside the Reissner-Nordstrom black hole horizons are 
nonvanishing, even though the Einstein scalar curvature vanishes in the interior 
of the charged metric.  Introducing a perfect fluid inside the Reissor-Nordstrom 
black hole, we have also explicitly evaluated the Einstein equations of motion inside the horizons 
to investigate the charge effects that the charge decreases the mass-energy density and the pressure 
of the fluid inside the black hole.

\vskip20pt
\noindent
{\bf ACKNOWLEDGMENTS}
\vskip10pt 

One of us (JC) would like to thank Prof. John K. Beem for helpful discussions and kind concerns.  
STH, JC and YJP would like to acknowledge financial support in part from Korea Science and 
Engineering Foundation (R01-2000-00015), (R01-2001-00003), and Ministry of Education, 
BK21 Project (D-1099), respectively.

\vskip18pt

${}^{1}$ W.G. Unruh, Phys. Rev. {\bf D14}, 870 (1976);
P.C.W. Davies, J. Phys. {\bf A8},\par\hskip0.3cm
 609 (1975).
\vskip5pt
${}^{2}$ S.W. Hawking, Comm. Math. Phys. {\bf 42}, 199 (1975).
\vskip5pt
${}^{3}$ E. Kasner, Am. J. Math. {\bf 43}, 130 (1921);
C. Fronsdal, Phys. Rev. {\bf 116}, 778\par\hskip0.3cm
 (1959).
\vskip5pt
${}^{4}$ J. Rosen, Rev. Mod. Phys. {\bf 37}, 204 (1965).
\vskip5pt
${}^{5}$ H. F. Goenner, {\em General Relativity and
Gravitation}, edited by A. Held (Plenum,\par\hskip0.3cm New York, 1980) 441.
\vskip5pt
${}^{6}$ H. Narnhofer, I. Peter  and W. Thirring, Int. J. Mod. Phys.
{\bf B10}, 1507 (1996).
\vskip5pt
${}^{7}$ S. Deser and O. Levin, Class. Quant. Grav. {\bf 14},
L163 (1997); Class. Quant.\par\hskip0.3cm Grav. {\bf 15}, L85 (1998); Phys. Rev. 
{\bf D59}, 0640004 (1999).
\vskip5pt
${}^{8}$ M. Banados, C. Teitelboim and J. Zanelli, Phys. Rev. Let.
{\bf 69}, 1849 (1992).
\vskip5pt
${}^{9}$  M. Banados, M. Henneaux, C. Teitelboim and J. Zanelli, Phys. Rev. {\bf D48},\par\hskip0.3cm
1506 (1993).
\vskip5pt
${}^{10}$ S. Carlip, Class. Quant. Grav. {\bf 12}, 2853 (1995).
\vskip5pt
${}^{11}$ D. Cangemi, M. Leblanc and R.B. Mann, Phys. Rev. {\bf D48}, 3606 (1993).
\vskip5pt
${}^{12}$ S.T. Hong, Y.W. Kim and Y.J. Park, Phys. Rev. {\bf D62}, 024024 (2000).
\vskip5pt
${}^{13}$ S.T. Hong, W.T. Kim, Y.W. Kim and Y.J. Park, Phys. Rev. {\bf D62}, 064021\par\hskip0.3cm (2000).
\vskip5pt
${}^{14}$ L. Andrianopoli, M. Derix, G.W. Gibbons, C. Herdeiro, A. Santambrogio and \par\hskip0.3cm
A. V. Proeyen, Class. Quant. Grav. {\bf 17}, 1875 (2000).
\vskip5pt
${}^{15}$  S. W. Hawking and H. S. Reall, Phys. Rev. {\bf D61}, 024014 (1999).
\vskip5pt
${}^{16}$ Y.W. Kim, Y.J. Park and K.S. Soh, Phys. Rev. {\bf D62}, 104020 (2000).
\vskip5pt
${}^{17}$ R.L. Bishop and B. O'Neill, Am. Math. Soc. {\bf 145}, 1 (1969).
\vskip5pt
${}^{18}$ R.L. Bishop and B. O'Neill, Trans. A.M.S. {\bf 145}, 1 (1969)
\vskip5pt
${}^{19}$ J.K. Beem, P. E. Ehrlich and K. Easley, {\it Global Lorentzian Geometry} 
(Marcel\par\hskip0.4cm Dekker Pure and Applied Mathematics, New York, 1996);
T.G. Powell, Ph.D\par\hskip0.4cm
Thesis, University of Missouri-Columbia (1982); K. Easley, Ph.D Thesis \par\hskip0.4cm
University of Missouri-Columbia (1980).
\vskip5pt
${}^{20}$ J.K. Beem and P. E. Ehrlich, Math. Proc. Camb. Phil. Soc. {\bf 85}, 161 (1979).
\vskip5pt
${}^{21}$ J. Choi, J. Math. Phys. {\bf 41}, 8163 (2000).
\vskip5pt
${}^{22}$ J.L. Flores and M. S\'{a}nchez, {\tt math.DG/9909075}; {\tt math.DG/0106174}. 
\vskip5pt
${}^{23}$ S.G. Harris, Class. Quant. Grav. {\bf 17}, 551 (2000).
\vskip5pt
${}^{24}$ C.W. Misner, K. S. Thorne and J. A. Wheeler, {\it Gravitation} (W. H. Freeman\par\hskip0.3cm
and Company, New York, 1970).
\vskip5pt
\end{document}